\documentclass[9pt,shortpaper,twoside,web]{IEEEtran}
\usepackage{cite}

\ifCLASSINFOpdf
  \usepackage[pdftex]{graphicx}
  \usepackage{epstopdf}
  \usepackage{subfigure}
  \graphicspath{{./}{./figures/}}
\else
\fi
\usepackage{xcolor}

\usepackage{amsmath}
\usepackage{amssymb}

\newcommand*{\qed}{\hfill\ensuremath{\square}}%
\newcommand*\diff{\mathop{}\!\mathrm{d}}%

\newtheorem{theorem}{Theorem}[section]
\newtheorem{corollary}{Corollary}[theorem]

\newtheorem{remark}[theorem]{Remark}
\newtheorem{definition}[theorem]{Definition}
\newtheorem{example}[theorem]{Example}

\allowdisplaybreaks

\begin{document}
%
\title{On the derivation of stability properties for time-delay systems without constraint on the time-derivative of the initial condition}
%
%
%

\author{Hugo~Lhachemi, Robert Shorten
\thanks{This publication has emanated from research supported in part by a research grant from Science Foundation Ireland (SFI) under grant number 16/RC/3872 and is co-funded under the European Regional Development Fund and by I-Form industry partners.
\newline\indent Hugo Lhachemi is with L2S, CentraleSup{\'e}lec, 91192 Gif-sur-Yvette, France (e-mail: hugo.lhachemi@centralesupelec.fr).
\newline\indent Robert Shorten is with the Dyson School of Design Engineering, Imperial College London, London, U.K, and also with the School of Electrical and Electronic Engineering, University College Dublin, Dublin, Ireland (e-mail: r.shorten@imperial.ac.uk).
}
}


%
%

\markboth{}%
{Lhachemi \MakeLowercase{\textit{et al.}}}
%



\maketitle

\begin{abstract}
Stability of retarded differential equations is closely related to the existence of Lyapunov-Krasovskii functionals. Even if a number of converse results have been reported regarding the existence of such functionals, there is a lack of constructive methods for their selection. For certain classes of time-delay systems for which such constructive methods are lacking, it was shown that Lyapunov-Krasovskii functionals that are also allowed to depend on the time-derivative of the state-trajectory are efficient tools for the study of the stability properties. However, in such an approach the initial condition needs to be assumed absolutely continuous with a square integrable weak derivative. In addition, the stability results hold for initial conditions that are evaluated based on the magnitude of both the initial condition and its time-derivative. The main objective of this paper is to show that, for certain classes of time-delay systems, the aforementioned stability results can actually be extended to initial conditions that are only assumed continuous and that are evaluated in uniform norm.
\end{abstract}


%
\IEEEpeerreviewmaketitle

\section{Introduction}\label{sec: Introduction}
%
%
%
%

The subject of network control systems is now a topic of great maturity in the control community \cite{hetel2017recent,liu2019survey}. Even though many results have been obtained for feedback loops operating over networks, some issues remain to be resolved. One such issue concerns feedback loops in which time-delays are time-varying, giving rise to certain types of delay differential equations. Such systems typically arise in networks with buffers and queues, which may or may not be filling, depending on levels of network congestion. For example, traditional communication networks, operating TCP flows, typically give rise to such effects. A basic question arising in the study of such networks concerns the effect of the regularity of initial conditions. As the systems may have been operating in delay free mode, before delays appear, initial conditions, and their effects, become a pertinent question. This question is even more pronounced in systems that may involve nonlinearity and switches, but also state resets~\cite{beker2001plant,hespanha2002switching}. Remarkably, this topic appears to have been only of secondary concern in the control community as very few results are reported~\cite{elsgoltz1973introduction,pepe2007problem}. The objective of this note is to derive new basic results on the effect of the regularity of initial conditions on available results on certain classes of retarded differential equations~\cite{elsgoltz1973introduction,kolmanovskii2012applied,hale2013introduction}. 

The stability of retarded differential equations is closely related to the existence of Lyapunov-Krasovskii functionals~\cite{gu2003stability}. This is why, during the past three decades, the use of Lyapunov-Krasovskii functionals has emerged as an efficient approach for the derivation of sufficient conditions to assess the stability properties of retarded differential equations~\cite{kolmanovskii2012applied,fridman2014introduction,fridman2014tutorial}. 
Even if converse results regarding the existence of such functionals have been reported for a large number of classes of time-delay systems, most of these results do not provide constructive methods for the selection of Lyapunov-Krasovskii functionals~\cite{gu2003stability,haidar2015converse}. For certain classes of time-delay systems for which such constructive methods are lacking, it was shown that Lyapunov-Krasovskii functionals that are also allowed to depend on the time-derivative of the state-trajectory are efficient tools for the study of the stability properties~\cite{fridman2001new,fridman2006new,lhachemi2019lmi}. In such a configuration the set of admissible initial conditions $W$ needs to be restricted to absolutely continuous functions $f : [-h_M,0] \rightarrow \mathbb{R}^n$ with square integrable weak derivative and is endowed with the norm $\Vert f \Vert_W = \sqrt{\Vert f(0) \Vert^2 + \int_{-h_M}^0 \Vert \dot{f}(s) \Vert^2 \,\mathrm{d}s}$. Consequently, the derived stability results are \emph{a priori} only valid for initial conditions $x_0 \in W$ and are generally stated for initial conditions evaluated by $\Vert x_0 \Vert_W$, i.e. in function of the magnitude of their time-derivative. 

The evaluation of the initial condition in $W$-norm and the restriction of the initial conditions to the set $W$ might be restrictive. First, in the context of time-delay systems, it is generally more natural to evaluate the initial conditions in uniform norm. However, $W$-norm and uniform norm are not equivalent over $W$. Indeed, while subsets of $W$ that are bounded in $W$-norm are also bounded in uniform norm, the converse is not true. This observation shows that the stability properties for initial conditions evaluated in uniform norm cannot be deduced from the ones for initial conditions evaluated in $W$-norm based on sole topological arguments. Second, the trajectories of time-delay systems are well-defined for a larger set of initial condition than $W$, namely the set of continuous functions. So, a natural question is whether the stability result can be 1) stated for initial conditions evaluated in uniform norm rather than in $W$-norm; 2) extended from initial conditions in $W$ to initial conditions that are assumed continuous. This paper gives a positive response to these two questions for a class of nonlinear systems. The derived result is obtained by leveraging a smoothing effect of the system trajectories for retarded differential equations.

It is worth noting that a similar issue, namely the choice of the norm in the evaluation of the state trajectory, was discussed in~\cite[Chap.~III Sec.~9]{elsgoltz1973introduction} in the context of neutral type systems. Specifically, it was shown under certain assumptions a form of equivalence between the evaluation of the state vector based on its only magnitude or also via the magnitude of its time-derivative. However, in both cases, the stability properties were stated for initial conditions evaluated based on both their magnitude and the magnitude of their time-derivative.

This paper is organized as follows. First, Section~\ref{sec: problem setting} provides the context and the main motivation of this work. Then the main results of this paper are presented in Section~\ref{sec: main result}. Finally, concluding remarks are provided in Section~\ref{sec: conclusion}.

\textbf{Notation} We endow $\mathbb{R}^n$ with the usual Euclidean norm $\Vert x \Vert = \sqrt{x^\top x}$. The set of continuous and continuously differentiable functions $f : I \rightarrow E$ from an interval $I \subset \mathbb{R}$ to a normed space $E$ is denoted by $\mathcal{C}^0(I;E)$ and $\mathcal{C}^1(I;E)$, respectively. For a given $h_M > 0$, the set $\mathcal{C}^0([-h_M,0];\mathbb{R}^n)$ is endowed with the uniform norm $\Vert f \Vert_\infty = \sup_{x \in [-h_M,0]} \Vert f(x) \Vert$. For any $t_f > 0$ and any continuous function $x : [-h_M , t_f ) \rightarrow \mathbb{R}^n$, we introduce for $0 \leq t < t_f$ the function $x_t = x(t+\cdot) \in \mathcal{C}^0([-h_M,0];\mathbb{R}^n)$. For a given constant $h_M > 0$, we introduce $W = W_{h_M}$ the set of absolutely continuous functions $f : [-h_M,0] \rightarrow \mathbb{R}^n$ with square integrable weak derivative. For any $a< b$, we recall that $g \in L^1(a,b)$ is the weak derivative of $f \in L^1(a,b)$ if $\int_a^b f(s) \varphi'(s) \,\mathrm{d}s = - \int_a^b g(s) \varphi(s) \,\mathrm{d}s$ for all smooth function $\varphi\in\mathcal{C}^\infty([a,b];\mathbb{R})$ with $\varphi(a)=\varphi(b)=0$. In this case, we note $\dot{f} \triangleq g$. We endow $W$ with the norm\footnote{This definition of the $W$-norm is extracted from~\cite[Chap.~4, Thm.~1.6, p~129]{kolmanovskii2012applied} in the context of the statement of the Lyapunov-Krasovskii theorem which is recalled in Theorem~\ref{thm: theorem of reference} of this paper. Note that certain authors prefer to adopt the following definition~\cite{fridman2014tutorial}: $\Vert f \Vert_W = \Vert f \Vert_\infty + \sqrt{\int_{-h_M}^0 \Vert \dot{f}(s) \Vert^2 \,\mathrm{d}s}$. However, such a choice has no impact on the resulting stability properties of the studied systems because standard computations show that the two aforementioned norms are equivalent in $W$.} $\Vert f \Vert_W = \sqrt{\Vert f(0) \Vert^2 + \int_{-h_M}^0 \Vert \dot{f}(s) \Vert^2 \,\mathrm{d}s}$. Finally, we introduce the following classes of comparison functions. A continuous function $\alpha : \mathbb{R}_+ \rightarrow \mathbb{R}_+$ is said to belong to the class $\mathcal{K}$ if it is strictly increasing with $\alpha(0)=0$. Moreover, $\alpha \in \mathcal{K}_\infty$ if $\alpha \in \mathcal{K}$ with $\lim_{s \rightarrow +\infty} \alpha(s) = +\infty$. A continuous function $\beta : \mathbb{R}_+ \times \mathbb{R}_+ \rightarrow \mathbb{R}_+$ is said to belong to the class $\mathcal{KL}$ if the mapping $s \mapsto \beta(s,t)$ is of class $\mathcal{K}$ for each fixed $t \geq 0$ and the mapping $t \mapsto \beta(s,t)$ is strictly decreasing with $\lim_{t \rightarrow +\infty} \beta(s,t) = 0$ for each fixed $s > 0$.

\section{Context and motivating example}\label{sec: problem setting}

\subsection{Definitions}

Let $f : A \times \mathcal{C}^0([-h_M,0];\mathbb{R}^n) \rightarrow \mathbb{R}^n$, with $A \subset \mathbb{R}^r$, be a continuous function such that $f(\cdot,0)=0$ and $f$ is locally Lipschitz continuous with respect to the second argument. Let $\mathcal{H}$ be any subset of $\mathcal{C}^0(\mathbb{R}_+; A)$. Let $X$ be a subspace of $\mathcal{C}^0([-h_M,0];\mathbb{R}^n)$ endowed with a norm denoted by $\Vert \cdot \Vert_X$. In the subsequent developments, $X$ will always be either $\mathcal{C}^0([-h_M,0];\mathbb{R}^n)$ endowed with the uniform norm or $W$ endowed with the $W$-norm. We consider the retarded differential equation:
\begin{subequations}\label{eq: RDE0}
\begin{align}
\dot{x}(t) & = f(\zeta(t),x_t) , & t \geq 0 \label{eq: RDE0 - equation} \\
x(\tau) & = x_0(\tau) , & \tau\in[-h_M,0] \label{eq: RDE0 - IC}
\end{align}
\end{subequations}
where $x_0 \in X$ and $\zeta \in \mathcal{H}$. We adopt the following definitions.

\begin{definition}[$(X,\mathcal{H})$-stability]\label{def: stab} System (\ref{eq: RDE0}) is said $(X,\mathcal{H})$-stable if there exist $\alpha \in \mathcal{K}$ and $R > 0$ such that for any $x_0 \in X$ with $\Vert x_0 \Vert_X < R$ and any $\zeta \in \mathcal{H}$, we have that the system trajectory $x$ is well defined on $\mathbb{R}_+$ with $\Vert x(t) \Vert \leq \alpha(\Vert x_0 \Vert_X)$ for all $t \geq 0$.
\end{definition}

\begin{definition}[$(X,\mathcal{H})$-local/global uniform asymptotic stability]\label{def: asymp stab} System (\ref{eq: RDE0}) is said $(X,\mathcal{H})$-locally uniformly asymptotically stable if there exists $\beta\in\mathcal{KL}$ and $R > 0$ such that for any $x_0 \in X$ with $\Vert x_0 \Vert_X < R$ and any $\zeta \in \mathcal{H}$ we have $\Vert x(t) \Vert \leq \beta(\Vert x_0 \Vert_X,t)$ for all $t \geq 0$. If one can select $R = + \infty$, system (\ref{eq: RDE0}) is said $(X,\mathcal{H})$-globally uniformly asymptotically stable.
\end{definition}

\begin{definition}[$(X,\mathcal{H})$-local/global exponential stability]\label{def: exp stab} System (\ref{eq: RDE0}) is said $(X,\mathcal{H})$-locally/globally exponentially stable if there exist constants $\kappa , C_0 > 0$ such that (\ref{eq: RDE0}) is $(X,\mathcal{H})$-locally/globally uniformly asymptotically stable with $\beta(s,t) = C_0 e^{-\kappa t} s$.
\end{definition}

\subsection{Motivating example}

In order to highlight the relevance of the study reported in this manuscript, we consider the below motivating example. Assume that we are interested in studying, for some constant $h_M > 0$, the stability properties of the following time-delay system:
\begin{subequations}\label{eq: sys1}
\begin{align}
\dot{x}(t) & = M x(t) + N x(t-h(t)) , & t \geq 0 \label{eq: sys1 - eq} \\
x(\tau) & = x_0(\tau) , & -h_M \leq \tau \leq 0  \label{eq: sys1 - IC}
\end{align}
\end{subequations}
where $M,N \in \mathbb{R}^{n \times n}$, $0 \leq h(t) \leq h_M$ is a continuous time-varying delay, and $x_0$ is the initial condition. Then one can resort to the following result\footnote{Note that this result has been improved later on in various directions~\cite{fridman2014introduction}. However, as the problematic studied in this paper is still relevant for these later developments, we motivate our study based on the result of the seminal work~\cite{fridman2001new} for ease and conciseness of the presentation.} borrowed from~\cite[Cor.~1]{fridman2001new} and \cite{fridman2014tutorial}.

\begin{theorem}\label{thm: stab sys1}
Let $M,N \in \mathbb{R}^{n \times n}$ and $h_M > 0$ be given. Assume that there exist $P_1 , Q \succ 0$ and $P_2,P_3 \in \mathbb{R}^{n \times n}$ such that 
\begin{equation}\label{eq: LMI}
\begin{bmatrix}
\Gamma^\top P_2 + P_2^\top \Gamma & P_1 -P_2^\top + \Gamma^\top P_3 & h_M P_2^\top N \\
P_1 - P_2 + P_3^\top \Gamma & - P_3 - P_3^\top + h_M Q & h_M P_3^\top N \\
h_M N^\top P_2 & h_M N^\top P_3 & - h_M Q
\end{bmatrix}
\prec 0 ,
\end{equation}
where $\Gamma = M+N$. Then, defining $\mathcal{H} = \{ h\in\mathcal{C}^0(\mathbb{R}_+;\mathbb{R}) \, : \, 0 \leq h(t) \leq h_M \}$, system (\ref{eq: sys1}) is $(W,\mathcal{H})$-globally uniformly asymptotically stable.
\end{theorem}

The reason of both the restriction of the initial conditions to the set $W$ and their evaluation in $W$-norm relies on the fact that, following~\cite{fridman2001new}, Theorem~\ref{thm: stab sys1} is proved by using the following Lyapunov-Krasovskii functional that depends on the time-derivative of the system trajectory: 
\begin{equation*}
V(t) = x(t)^\top P_1 x(t) + \int_{-h_M}^0 \int_{t+\theta}^{t} \dot{x}(s)^\top Q \dot{x}(s) \,\mathrm{d}s\,\mathrm{d}\theta .   
\end{equation*}
For times $0 \leq t < h_M$, the integral term of the above Lyapunov-Krasovskii functional involves the time derivative of the initial condition $x_0$. This is the main motivation for considering initial conditions $x_0 \in W$, ensuring that $V$ is well defined for all $t \geq 0$. Then, it was shown in~\cite{fridman2001new} under the assumptions of Theorem~\ref{thm: stab sys1} that there exists $\gamma > 0$, independent of $x_0$ and $h$, such that $\dot{V}(t) \leq - \gamma \Vert x(t) \Vert^2$. Finally, the asymptotic stability result stated in Theorem~\ref{thm: stab sys1} follows from the existence of constants $\alpha,\beta>0$, independent of $x_0$ and $h$, such that $\alpha \Vert x (t) \Vert^2 \leq V(t) \leq \beta \Vert x (t+\cdot) \Vert_W^2$. The latter estimate is the reason why the initial condition is evaluated in $W$-norm in the statement of Theorem~\ref{thm: stab sys1}.

In this context, the question addressed in this paper is whether, for certain classes of systems (which will include (\ref{eq: sys1})), a $(W,\mathcal{H})$-stability property implies a $(\mathcal{C}^0([-h_M,0];\mathbb{R}^n),\mathcal{H})$-stability property. The motivation behind this question is twofold. First, while the initial stability properties hold for initial conditions in $W$, the trajectories of (\ref{eq: sys1}) are actually well defined for the strictly larger\footnote{For instance, the function $f(t) = \sqrt{\vert t \vert}$ is absolutely continuous over $[-h_M,0]$ but its derivative is not square integrable over $[-h_M,0]$. There also exist functions that are continuous but not absolutely continuous. Two classical examples are the Cantor function and the function $f(t) = t \sin(1/t)$ for $t \neq 0$ and $f(0) = 0$ in any compact interval containing $0$.} set of initial conditions $\mathcal{C}^0([-h_M,0];\mathbb{R}^n)$. Second, $W$-norm and uniform norm are not equivalent on $W$. Indeed, while subsets of $W$ that are bounded in $W$-norm are also bounded in uniform norm, because we have $\Vert f \Vert_\infty \leq (1+\sqrt{h_M}) \Vert f \Vert_W$ for all $f \in W$, the converse is not true. A classical example of this is provided below. This shows that the derivation of  a $(\mathcal{C}^0([-h_M,0];\mathbb{R}^n),\mathcal{H})$-stability property from a pre-existing $(W,\mathcal{H})$-stability property cannot be deduced based on sole topological arguments.

\begin{example}
Consider the function $\varphi:[0,1]\rightarrow\mathbb{R}$ defined by $\varphi(t) = 2t$ if $0 \leq t \leq 1/2$ and $\varphi(t) = 2(1-t)$ if $1/2 \leq t \leq 1$. Then define $\phi(t) = \varphi(t-\lfloor t \rfloor)$ for any $t \in \mathbb{R}$, where we recall that $\lfloor \cdot \rfloor$ is the floor function, i.e. $\lfloor t \rfloor$ is the unique integer such that $\lfloor t \rfloor \leq t < \lfloor t \rfloor +1$. Finally, for a given unit vector $x \in \mathbb{R}^n$, consider for any $m \geq 1$ the function $f_m(t) = \phi(m t / h_M) x$ for $-h_M \leq t \leq 0$. Then one can see that $f_m \in W$, $\Vert f_m \Vert_\infty = 1$, and $\Vert f_m \Vert_W = 2 m / \sqrt{h_M} \rightarrow +\infty$ when $m \rightarrow + \infty$. This shows that there cannot exist a constant $C > 0$ such that $\Vert f \Vert_W \leq C \Vert f \Vert_\infty$ for all $f \in W$. This implies that $W$-norm and uniform norm are not equivalent in $W$.
\end{example}

\section{Main results}\label{sec: main result}

\subsection{Sufficient conditions ensuring the equivalence of $(W,\mathcal{H})$ and $(\mathcal{C}^0([-h_m,0];\mathbb{R}^n),\mathcal{H})$ stability properties}

We start with the introduction of the following definition.

\begin{definition}\label{def: forward invariance}
We say that $\mathcal{H} \subset \mathcal{C}^0(\mathbb{R}_+;\mathbb{R}^r)$ is forward invariant if for any $\zeta \in \mathcal{H}$ and any $\alpha > 0$ we have $\zeta(\cdot + \alpha) \in \mathcal{H}$.
\end{definition}

The set $\mathcal{H} \subset \mathcal{C}^0(\mathbb{R}_+;\mathbb{R}^r)$ is used as the set of admissible functions $\zeta$ for the time-delay system (\ref{eq: RDE0}). Considering systems exhibiting time-delays, Definition~\ref{def: forward invariance} is satisfied by most of the relevant applications. In particular, the subset of $\mathcal{C}^0(\mathbb{R}_+;\mathbb{R}^r)$ composed of functions $h = (h_1,\ldots,h_r)$ such that $0 \leq h_i(t) \leq h_M$ is forward invariant.

We can now introduce the following main result.

\begin{theorem}\label{thm: transfer theorem}
Let $f : A \times \mathcal{C}^0([-h_M,0];\mathbb{R}^n) \rightarrow \mathbb{R}^n$, with $A \subset \mathbb{R}^r$, be a continuous function such that $f(\cdot,0)=0$ and $f$ is locally Lipschitz continuous with respect to the second argument. Let $\mathcal{H}$ be any subset of $\mathcal{C}^0(\mathbb{R}_+; \mathbb{R}^r)$ that is forward invariant and that satisfies $\zeta(t) \in A$ for all $\zeta \in \mathcal{H}$ and all $t \geq 0$. We consider the retarded differential equation:
\begin{subequations}\label{eq: RDEtransfer}
\begin{align}
\dot{x}(t) & = f(\zeta(t),x_t) , & t \geq 0 \label{eq: RDEtransfer - equation} \\
x(\tau) & = x_0(\tau) , & \tau\in[-h_M,0] \label{eq: RDEtransfer - IC}
\end{align}
\end{subequations}
where $x_0 \in \mathcal{C}^0([-h_M,0];\mathbb{R}^n)$ and $\zeta \in \mathcal{H}$. For a given $H > 0$, assume that the restriction of $f$ to $A \times \{ \phi \in \mathcal{C}^0([-h_M,0];\mathbb{R}^n) \,:\, \Vert \phi \Vert_\infty \leq H \}$ is Lipschitz continuous with respect to the second argument, uniformly with respect to the first argument. Then the following hold for system (\ref{eq: RDEtransfer}):
\begin{itemize}
\item[(i)] $(W,\mathcal{H})$-stability $\Leftrightarrow$ $(\mathcal{C}^0([-h_M,0];\mathbb{R}^n),\mathcal{H})$-stability;

\item[(ii)] $(W,\mathcal{H})$-local uniform asymptotic stability $\Leftrightarrow$ $(\mathcal{C}^0([-h_M,0];\mathbb{R}^n),\mathcal{H})$-local uniform asymptotic stability;

\item[(iii)] $(W,\mathcal{H})$-local exponential stability $\Leftrightarrow$ $(\mathcal{C}^0([-h_m,0];\mathbb{R}^n),\mathcal{H})$-local exponential stability.
\end{itemize}
Moreover, if $f$ is globally Lipschitz continuous with respect to the second argument, uniformly with respect to the first argument, then (ii) and (iii) still hold when replacing ``locally'' by ``globally''.
\end{theorem}

\textbf{Proof.}
The ``$\Leftarrow$'' parts are a trivial consequence of $W \subset \mathcal{C}^0([-h_m,0];\mathbb{R}^n)$ and that $\Vert f \Vert_\infty \leq (1+\sqrt{h_M}) \Vert f \Vert_W$ for all $f \in W$. Thus we focus on the ``$\Rightarrow$'' parts.

Let $L>0$ be such that $\Vert f(\xi,\phi) \Vert \leq L \Vert \phi \Vert_\infty$ for all $\xi \in A$ and all $\phi \in \mathcal{C}^0([-h_M,0];\mathbb{R}^n)$ with $\Vert \phi \Vert_\infty \leq H$. Let $x_0 \in \mathcal{C}^0([-h_M,0];\mathbb{R}^n)$ and $\zeta \in \mathcal{H}$ be given. Noting that $(t,\phi) \mapsto f(\zeta(t),\phi)$ is continuous and locally Lipschitz continuous with respect to the second argument, we can introduce $x : [-h_M,t_f) \rightarrow \mathbb{R}^n$ the maximal solution of (\ref{eq: RDEtransfer}) where $0 < t_f \leq + \infty$. 

In a first time we show that $x_0 \in \mathcal{C}^0([-h_M,0];\mathbb{R}^n)$ with $\Vert x_0 \Vert_\infty < H / (1+e^{L h_M})$ implies $t_f > h_M$ and $\Vert x_t \Vert_\infty < H$ for all $0 \leq t \leq h_M$. There are two cases. In the first case we have $\Vert x_t \Vert_\infty < H$ for all $t \in [0,t_f)$. As $\Vert f(\zeta(t),\phi) \Vert \leq L H$ for all $t \geq 0$ and all $\phi \in \mathcal{C}^0([-h_M,0];\mathbb{R}^n)$ with $\Vert \phi \Vert_\infty \leq H$, we obtain (see e.g. \cite[Chap.~2, Thm.~3.2]{hale2013introduction}) that $t_f = + \infty$. In the second case, there exists $0 \leq t_1 < t_f$ such that $\Vert x_{t_1} \Vert_\infty \geq H$. Since $\Vert x_{0} \Vert_\infty < H$ and $t \mapsto \Vert x_t \Vert_\infty $ is continuous, we can select $t_1 > 0$ such that $\Vert x_t \Vert_\infty < H$ for all $0 \leq t < t_1$ and $\Vert x_{t_1} \Vert_\infty = H$. Then we have, for all $0 \leq t \leq t_1$,
\begin{align*}
\Vert x(t) \Vert
& \leq \Vert x(0) \Vert + \int_0^t \Vert \dot{x}(\tau) \Vert \,\diff\tau \\
& \leq \Vert x_0 \Vert_\infty + \int_0^t \Vert f(\zeta(\tau),x_\tau) \Vert \,\diff\tau \\
& \leq \Vert x_0 \Vert_\infty + L \int_0^t \Vert x_\tau \Vert_\infty \,\diff\tau ,
\end{align*}
and thus
\begin{align*}
\Vert x_t \Vert_\infty 
& \leq \max\left( \sup\limits_{\tau \in [-h_M,0]} \Vert x(\tau) \Vert , \sup\limits_{\tau \in [\max(t-h_M,0),t]} \Vert x(\tau) \Vert \right) \\
& \leq \Vert x_0 \Vert_\infty + L \int_0^t \Vert x_\tau \Vert_\infty \,\diff\tau .
\end{align*}
The use of Gr{\"o}nwall's inequality implies
\begin{equation*}
\Vert x_t \Vert_\infty
\leq (1+e^{L t}) \Vert x_0 \Vert_\infty 
\leq (1+e^{L t_1}) \Vert x_0 \Vert_\infty 
\end{equation*}
for all $0 \leq t \leq t_1$. Recalling that $\Vert x_0 \Vert_\infty < H / (1+e^{L h_M})$, the use of $\Vert x_{t_1} \Vert_\infty = H$ gives $h_M < t_1 < t_f$. In both cases, $\Vert x_0 \Vert_\infty < H / (1+e^{L h_M})$ implies $t_f > h_M$ and the replication of the above reasoning shows that $\Vert x_t \Vert_\infty \leq (1+e^{L h_M}) \Vert x_0 \Vert_\infty < H$ for all $0 \leq t \leq h_M$. 

In a second time we show that, for any $\gamma > 0$, $\Vert x_0 \Vert_\infty < \dfrac{1}{1+e^{L h_M}} \min\left( H , \frac{\gamma}{\sqrt{1+h_M L^2}} \right)$ implies $\Vert x_{h_M} \Vert_W < \gamma$. Indeed, as the right hand side of (\ref{eq: RDEtransfer - equation}) is continuous, we have for $\Vert x_0 \Vert_\infty < H / (1+e^{L h_M})$ that $x_{h_M} \in W$ and  
\begin{align}
\Vert x_{h_M} \Vert_W
& = \sqrt{ \Vert x(h_M) \Vert^2 + \int_0^{h_M} \Vert \dot{x}(\tau) \Vert^2 \,\diff\tau} \nonumber \\
& \leq \sqrt{ \Vert x_{h_M} \Vert_\infty^2 + L^2 \int_0^{h_M} \Vert x_\tau \Vert_\infty^2 \,\diff\tau} \nonumber \\
& \leq (1+e^{L h_M}) \sqrt{ 1 + h_M L^2 } \Vert x_0 \Vert_\infty , \label{eq: estimate W norm x_hM}
\end{align}
which provides the claimed conclusion.

Now, for any initial condition $x_0 \in \mathcal{C}^0([-h_M,0];\mathbb{R}^n)$ with $\Vert x_0 \Vert_\infty < H / (1+e^{L h_M})$, which implies $t_f > h_M$, we can introduce $\tilde{x}(t) = x(t+h_M)$ for $t\in[-h_M,t_f-h_M)$, $\tilde{\zeta}(t) = \zeta(t+h_M)$ for $t \geq 0$, and $\tilde{x}_0=x_{h_M} \in W$. As
\begin{align*}
\dot{\tilde{x}}(t)
= \dot{x}(t+h_M)
= f(\zeta(t+h_M),x_{t+h_M})
= f(\tilde{\zeta}(t),\tilde{x}_{t})
\end{align*}
for all $0 \leq t < t_f - h_M$ with the initial condition 
\begin{equation*}
\tilde{x}(\tau) = x(h_M+\tau) = x_{h_M}(\tau) = \tilde{x}_0(\tau)
\end{equation*}
for $\tau\in[-h_M,0]$, then $\tilde{x}$ is the maximal solution\footnote{Otherwise we could extend $x$, which we cannot since $x$ is assumed to be a maximal solution.} of (\ref{eq: RDEtransfer}) associated with the initial condition $\tilde{x}_0 \in W$ and $\tilde{\zeta} \in \mathcal{H}$ because $\mathcal{H}$ is forward invariant.

Assume now that (\ref{eq: RDEtransfer}) is $(W,\mathcal{H})$-locally uniformly asymptotically stable, i.e., there exists $\beta\in\mathcal{KL}$ and $R > 0$ such that for any $x_0 \in W$ with $\Vert x_0 \Vert_W < R$ and any $\zeta \in \mathcal{H}$, the associated system trajectory satisfies $\Vert x(t) \Vert \leq \beta(\Vert x_0 \Vert_W,t)$ for all $t \geq 0$. Then, for any initial condition $x_0 \in \mathcal{C}^0([-h_M,0];\mathbb{R}^n)$ such that $\Vert x_0 \Vert_\infty < \tilde{R} \triangleq \dfrac{1}{1+e^{L h_M}} \min\left( H , \frac{R}{\sqrt{1+h_M L^2}} \right)$, we have $\Vert \tilde{x}_0 \Vert_W = \Vert x_{h_M} \Vert_W < R$. Hence we infer that $t_f = + \infty$ and, for all $t \geq 0$, $\Vert \tilde{x} (t) \Vert \leq \beta ( \Vert \tilde{x}_0 \Vert_W , t )$. Based on (\ref{eq: estimate W norm x_hM}) we obtain that
\begin{equation*}
\Vert x(t) \Vert \leq \beta\left( (1+e^{L h_M}) \sqrt{ 1 + h_M L^2 } \Vert x_0 \Vert_\infty , t-h_M \right)
\end{equation*}
for all $t \geq h_M$ while we have already shown that $\Vert x(t) \Vert \leq \Vert x_t \Vert_\infty \leq (1+e^{L h_M}) \Vert x_0 \Vert_\infty$ for all $0 \leq t \leq h_M$. This shows the existence of $\tilde{\beta}\in\mathcal{KL}$, which only depends on $\beta$, $h_M$, and $L$, such that for any $x_0 \in \mathcal{C}^0([-h_M,0];\mathbb{R}^n)$ with $\Vert x_0 \Vert_\infty < \tilde{R}$ and any $\zeta \in \mathcal{H}$, the trajectory $x$ of (\ref{eq: RDEtransfer}) satisfies $\Vert x(t) \Vert \leq \tilde{\beta}(\Vert x_0 \Vert_\infty , t)$. This shows (ii). Points (i) and (iii) follow by similar arguments.

Finally, if one assumes that $f$ is globally Lipschitz continuous with respect to the second argument, uniformly with respect to the second argument, the above results hold with $H = +\infty$. In that case, this shows that (ii) and (iii) still hold when replacing ``locally'' by ``globally''.
\qed

Consider now the following class of linear time-delay systems:
\begin{subequations}\label{eq: lin RDE}
\begin{align}
\dot{x}(t) & = A x(t) + \sum\limits_{i = 1}^{p} B_i x(t-h_i(t)) \label{eq: lin RDE - equation} \\
& \phantom{=}\; + \sum\limits_{i = p+1}^{p+q} \int_{-h_i(t)}^{0} B_i(\theta) x(t+\theta) \,\mathrm{d}\theta , & t \geq 0 \nonumber \\
x(\tau) & = x_0(\tau) , & \tau\in [-h_M,0] \label{eq: lin RDE - IC}
\end{align}
\end{subequations}
with $h_M > 0$, $A,B_1,\ldots,B_p\in\mathbb{R}^{n \times n}$, $B_{p+1},\ldots,B_{p+q} : [-h_M,0] \rightarrow \mathbb{R}^{n \times n}$ integrable, and $h_1,\ldots,h_{p+q} \in \mathcal{C}^0(\mathbb{R}_+;\mathbb{R})$ with $0 \leq h_i(t) \leq h_M$. In this setting, $h_M > 0$ stands for a known upper-bound of the delays $h_i$. In this configuration, the following result is an immediate corollary of Theorem~\ref{thm: transfer theorem}.

\begin{corollary}\label{thm: main result}
Let $h_M > 0$, $A,B_1,\ldots,B_p\in\mathbb{R}^{n \times n}$, $B_{p+1},\ldots,B_{p+q} : [-h_M,0] \rightarrow \mathbb{R}^{n \times n}$ integrable, and $\mathcal{H} \subset \mathcal{C}^0(\mathbb{R}_+;\mathbb{R}^{p+q})$ that is forward invariant and such that any $h = (h_1,\ldots,h_{p+q}) \in \mathcal{H}$ satisfies $0 \leq h_i (t) \leq h_M$ for $1 \leq i \leq p+q$ and $t \geq 0$. Then, the followings hold for system (\ref{eq: lin RDE}):
\begin{itemize}
\item[(i)] $(W,\mathcal{H})$-stability $\Leftrightarrow$ $(\mathcal{C}^0([-h_m,0];\mathbb{R}^n),\mathcal{H})$-stability;

\item[(ii)] $(W,\mathcal{H})$-global uniform asymptotic stability $\Leftrightarrow$ $(\mathcal{C}^0([-h_m,0];\mathbb{R}^n),\mathcal{H})$-global uniform asymptotic stability;

\item[(iii)] $(W,\mathcal{H})$-global exponential stability $\Leftrightarrow$ $(\mathcal{C}^0([-h_m,0];\mathbb{R}^n),\mathcal{H})$-global exponential stability.
\end{itemize} 
\end{corollary}

\begin{remark}
For linear time-invariant systems it is well known that asymptotic stability implies exponential stability. In the context of linear time-delay systems with time-varying delays, the situation appears to be much more complex~\cite[p.~136]{fridman2014introduction}. Nevertheless, it was shown in~\cite{liu2013relationships} that such an implication holds for a class of positive linear systems with bounded time-varying delays.
\end{remark}

\begin{example}
The application of Corollary~\ref{thm: main result} to the setting of Theorem~\ref{thm: stab sys1} immediately yields the following result. Under the assumptions of Theorem~\ref{thm: stab sys1}, system (\ref{eq: sys1}) is $(\mathcal{C}^0([-h_M,0];\mathbb{R}^n),\mathcal{H})$-globally uniformly asymptotically stable.
\end{example}

\subsection{A Lyapunov-Krasovskii-based approach for the stability assessment of nonlinear retarded differential equations}\label{subsec: main thm nonlinear}

In this subsection we fix a constant $h_M > 0$ and we introduce $Q_H = \{ \phi \in W \,:\, \Vert \phi \Vert_W \leq H \}$. For a function $V:[a,b)\rightarrow\mathbb{R}$, the Dini upper right derivative is defined as
\begin{equation*}
D^+V(t) = \limsup\limits_{h \rightarrow 0^+} \dfrac{V(t+h)-V(t)}{h} 
\end{equation*}
for any $t \in [a,b)$.

One of the main results (e.g. used in~\cite{fridman2001new} to complete the proof of Theorem~\ref{thm: stab sys1}) for the study of the stability properties of retarded differential equations by means of Lyapunov-Krasovskii functionals depending on the time derivative of the system trajectories is given by the following theorem borrowed from~\cite[Chap.~4, Thm.~1.6]{kolmanovskii2012applied}.

\begin{theorem}\label{thm: theorem of reference}
Let $f : \mathbb{R}_+ \times \mathcal{C}^0([-h_M,0];\mathbb{R}^n) \rightarrow \mathbb{R}^n$ be a continuous function such that $f(\cdot,0)=0$ and $f$ is locally Lipschitz continuous with respect to the second argument. We consider the retarded differential equation:
\begin{subequations}\label{eq: RDE}
\begin{align}
\dot{x}(t) & = f(t,x_t) , & t \geq 0 \label{eq: RDE - equation} \\
x(\tau) & = x_0(\tau) , & \tau\in[-h_M,0] \label{eq: RDE - IC}
\end{align}
\end{subequations}
with initial condition $x_0 \in W$. For a given $H > 0$, assume that the restriction of $f$ to $\mathbb{R}_+ \times \{ \phi \in \mathcal{C}^0([-h_M,0];\mathbb{R}^n) \,:\, \Vert \phi \Vert_\infty \leq H \}$ is Lipschitz continuous with respect to the second argument, uniformly with respect to the first argument, and with Lipschitz constant denoted by $L > 0$. Let there exist a continuous functional $V : \mathbb{R}_+ \times W \rightarrow\mathbb{R}$ and $\omega_1 , \omega_2 \in \mathcal{K}$ such that
\begin{equation}\label{eq: theorem of reference - assumed inequality 1}
\omega_1(\Vert \phi(0) \Vert)
\leq V(t,\phi)
\leq \omega_2(\Vert \phi \Vert_W)
\end{equation}  
for all $t \geq 0$ and $\phi \in Q_H$. Introducing $V(t) = V(t,x_t)$ with $x$ the trajectory of (\ref{eq: RDE}), assume that 
\begin{equation}\label{eq: theorem of reference - assumed inequality 2}
D^+V \leq 0
\end{equation}
as soon as $x_t \in Q_H$. Then the origin of (\ref{eq: RDE}) is locally stable in the sense that there exist $\alpha\in\mathcal{K}$ and $R > 0$ such that for any $x_0 \in W$ with $\Vert x_0 \Vert_W < R$, the system trajectory $x$ of (\ref{eq: RDE}) is well defined for all $t \geq 0$ and $\Vert x(t) \Vert \leq \alpha(\Vert x_0 \Vert_W)$ for all $t \geq 0$. Moreover, if there exists $\omega_3 \in \mathcal{K}$ such that 
\begin{equation}\label{eq: theorem of reference - assumed inequality 3}
D^+ V(t) \leq - \omega_3(\Vert x(t) \Vert)
\end{equation}
as soon as $x_t \in Q_H$,  then the origin of (\ref{eq: RDE}) is locally uniformly asymptotically stable in the sense that there exist $\beta\in\mathcal{KL}$ and $R > 0$ such that for any $x_0 \in W$ with $\Vert x_0 \Vert_W < R$, the system trajectory $x$ of (\ref{eq: RDE}) is well defined for all $t \geq 0$ and $\Vert x(t) \Vert \leq \beta(\Vert x_0\Vert_W , t )$ for all $t \geq 0$. Furthermore, if the above assumptions are satisfied with $\omega_1,\omega_2 \in \mathcal{K}_\infty$ and for $H = +\infty$, then the uniform asymptotic stability result holds globally ($R = +\infty$).
\end{theorem}

\begin{remark}
In certain references, see e.g.~\cite[Chap.~4, Thm.~1.6]{kolmanovskii2012applied}, the dependency of the Lyapunov-Krasovskii functional $V$ on the time-derivative $\dot{x}_t$ of the system trajectory $x_t$ is explicitly highlighted by using the notation $V(t) = V(t,x_t,\dot{x}_t)$. In the statement of Theorem~\ref{thm: theorem of reference}, such a dependency is captured by the fact that the second argument of $V$ belongs to $W$ endowed with the $W$-norm.
\end{remark}

\begin{remark}
The version of Theorem~\ref{thm: theorem of reference} reported in~\cite[Chap.~4, Thm.~1.6]{kolmanovskii2012applied} is stated in terms of an asymptotic stability property. Nevertheless, the uniform asymptotic stability property also follows from a straightforward adaptation of the arguments reported in~\cite[Chap.~5, proof of Thm.~2.1]{hale2013introduction} in the context of Lyapunov-Krasovskii functional that only depend on the state trajectory (but not on its time derivative). Moreover, even if the developments reported in these references  are concerned with local asymptotic stability properties, their proofs readability extend to global asymptotic stability properties as soon as the vector field $f$ is assumed to be globally Lipschitz continuous with respect to the second argument. Finally, it is worth noting that the statement of the above stability results in terms of comparison functions (instead of quantifiers-based definitions of the stability properties) follows from classical results; see e.g.~\cite[proof of Lem.~4.5]{khalil2002nonlinear}.
\end{remark}

\begin{remark}\label{rem: thm of reference - independance of dela versus f and V}
Let $(f_j)_{j \in J}$ be a family (either finite of infinite, either countable or uncountable) of vector fields all satisfying the assumptions of Theorem~\ref{thm: theorem of reference} with possibly different Lyapunov-Krasovskii functionals $(V_j)_{j \in J}$ but with the same functions $\omega_1,\omega_2 \in \mathcal{K}$ and the same constants $h_M,H,L > 0$. It follows from the proof of Theorem~\ref{thm: theorem of reference} (see~\cite[Chap.~4, Thm.~1.6]{kolmanovskii2012applied} and~\cite[Chap.~5, Thm.~2.1]{hale2013introduction}) that the stability property holds uniformly w.r.t. $j \in J$. In other words, there exist $\alpha\in\mathcal{K}$ and $R > 0$ such that for all $j \in J$ and all $x_0 \in W$ with $\Vert x_0 \Vert_W < R$, the system trajectory $x_j$ of (\ref{eq: RDE}), when replacing $f$ by $f_j$, is well defined for all $t \geq 0$ and $\Vert x_j(t) \Vert \leq \alpha(\Vert x_0 \Vert_W)$ for all $t \geq 0$. Moreover, if there exists a function $\omega_3 \in \mathcal{K}$ such that, when replacing $V$ by $V_j$, (\ref{eq: theorem of reference - assumed inequality 3}) holds for all $j \in J$, then there exist $\beta\in\mathcal{KL}$ and $R > 0$ such that for all $j \in J$ and all $x_0 \in W$ with $\Vert x_0 \Vert_W < R$, the system trajectory $x_j$ of (\ref{eq: RDE}), when replacing $f$ by $f_j$, is well defined for all $t \geq 0$ and $\Vert x_j(t) \Vert \leq \beta(\Vert x_0\Vert_W , t )$ for all $t \geq 0$. Moreover, if one can further select $\omega_1,\omega_2 \in \mathcal{K}_\infty$ and $H = +\infty$, then this uniform asymptotic stability result holds globally ($R = +\infty$). This remark plays a key role in the proof of the next theorem.
\end{remark}

We now introduce an alternative version of Theorem~\ref{thm: theorem of reference} that extends the stability properties from initial conditions $x_0 \in W$ to $x_0 \in \mathcal{C}^0([-h_M,0];\mathbb{R}^n)$ and with initial conditions evaluated in uniform norm instead of the $W$-norm.

\begin{theorem}\label{thm: final theorem}
Let $f : A \times \mathcal{C}^0([-h_M,0];\mathbb{R}^n) \rightarrow \mathbb{R}^n$, with $A \subset \mathbb{R}^r$, be a continuous function such that $f(\cdot,0)=0$ and $f$ is locally Lipschitz continuous with respect to the second argument. Let $\mathcal{H}$ be any subset of $\mathcal{C}^0(\mathbb{R}_+; \mathbb{R}^r)$ that is forward invariant and that satisfies $\zeta(t) \in A$ for all $\zeta \in \mathcal{H}$ and all $t \geq 0$. We consider the retarded differential equation:
\begin{subequations}\label{eq: RDE2}
\begin{align}
\dot{x}(t) & = f(\zeta(t),x_t) , & t \geq 0 \label{eq: RDE2 - equation} \\
x(\tau) & = x_0(\tau) , & \tau\in[-h_M,0] \label{eq: RDE2 - IC}
\end{align}
\end{subequations}
where $x_0 \in \mathcal{C}^0([-h_M,0];\mathbb{R}^n)$ and $\zeta \in \mathcal{H}$. For a given $H > 0$, assume that the restriction of $f$ to $A \times \{ \phi \in \mathcal{C}^0([-h_M,0];\mathbb{R}^n) \,:\, \Vert \phi \Vert_\infty \leq H \}$ is Lipschitz continuous with respect to the second argument, uniformly with respect to the first argument. Let there exist a continuous functional $V : A \times W \rightarrow\mathbb{R}$ and $\omega_1 , \omega_2 \in \mathcal{K}$ such that
\begin{equation}\label{eq: main thm 2 - assumed inequality 1}
\omega_1(\Vert \phi(0) \Vert)
\leq V(\xi,\phi)
\leq \omega_2(\Vert \phi \Vert_W)
\end{equation}  
for all $\xi \in A$ and $\phi \in Q_H$. For any $x_0 \in W$ and $\zeta \in \mathcal{H}$, introducing $V_\zeta(t) = V(\zeta(t),x_t)$ with $x$ the trajectory of (\ref{eq: RDE2}), assume that 
\begin{equation}\label{eq: main thm 2 - assumed inequality 2}
D^+V_\zeta \leq 0
\end{equation}
as soon as $x_t \in Q_H$. Then (\ref{eq: RDE2}) is $(\mathcal{C}^0([-h_M,0];\mathbb{R}^n),\mathcal{H})$-stable. Moreover, if there exists $\omega_3 \in \mathcal{K}$ such that 
\begin{equation}\label{eq: main thm 2 - assumed inequality 3}
D^+ V_\zeta(t) \leq - \omega_3(\Vert x(t) \Vert)
\end{equation}
as soon as $x_t \in Q_H$, then (\ref{eq: RDE2}) is $(\mathcal{C}^0([-h_M,0];\mathbb{R}^n),\mathcal{H})$-locally uniformly asymptotically stable. Furthermore, if the above assumptions are satisfied with $\omega_1,\omega_2 \in \mathcal{K}_\infty$ and for $H = +\infty$, then the uniform asymptotic stability result holds globally.
\end{theorem}

\begin{remark}
It is worth noting that the functional $V$, as well as the assumed inequalities (\ref{eq: main thm 2 - assumed inequality 1}-\ref{eq: main thm 2 - assumed inequality 3}), are well-defined only for system trajectories associated with initial conditions $x_0 \in W$ and evaluated in $W$-norm. In this context, the main interest of Theorem~\ref{thm: final theorem} relies on the fact that the $(\mathcal{C}^0([-h_M,0];\mathbb{R}^n),\mathcal{H})$-stability property of the system is deduced from the study of the system trajectories associated with initial conditions exhibiting extra regularity assumption, namely $x_0 \in W$, and evaluated in $W$-norm. Note that the function $\zeta$ is only assumed continuous. Thus, in general, $t \mapsto V(t)$ is only continuous but not differentiable. This is not an issue since we consider in the statement of Theorem~\ref{thm: final theorem} the concept of Dini upper right derivative, which implies that the quantity $D^+V(t)$ is always well defined for all $t \geq 0$.
\end{remark}

\begin{remark}
Note that even if system (\ref{eq: RDE2}) can be cast into the form (\ref{eq: RDE}) for a given $\zeta \in \mathcal{H}$, the derived proof of Theorem~\ref{thm: final theorem} explicitly requires the structure (\ref{eq: RDE2}) and cannot be applied, in general, to (\ref{eq: RDE}). This is because the proof essentially relies on the application of Theorem~\ref{thm: transfer theorem} which explicitly uses the foward invariance property of the set $\mathcal{H}$. Such a procedure cannot be applied to system (\ref{eq: RDE}) because, in general, $\dot{\tilde{x}}(t)  = \dot{x}(t+h_M) = f(t+h_M,\tilde{x}_{t}) \neq f(t,\tilde{x}_{t})$. Then, we have in general that $\tilde{x}$ is not a solution of (\ref{eq: RDE}) associated with the initial condition $\tilde{x}_0$. This is the main reason why the proof reported hereafter cannot be applied to (\ref{eq: RDE}).
\end{remark}

\textbf{Proof.}
For any $\zeta \in \mathcal{H}$, we define $f_\zeta(t,\phi) = f(\zeta(t),\phi)$ for all $t \geq 0$ and $\phi \in \mathcal{C}^0([-h_M,0];\mathbb{R}^n)$. Hence the vector field $f_\zeta$ satisfies the assumptions of Theorem~\ref{thm: theorem of reference}. Applying Remark~\ref{rem: thm of reference - independance of dela versus f and V} to the family of vector fields $(f_\zeta)_{\zeta \in \mathcal{H}}$, we infer that (\ref{eq: RDE2}) is $(W,\mathcal{H})$-stable. Since $\mathcal{H}$ is assumed forward invariant, we conclude from Theorem~\ref{thm: transfer theorem} that (\ref{eq: RDE2}) is $(\mathcal{C}^0([-h_M,0];\mathbb{R}^n),\mathcal{H})$-stable. The same reasoning applies for the uniform asymptotic stability property.
\qed

\begin{remark}
Under the regularity assumptions of Theorem~\ref{thm: final theorem}, sufficient conditions ensuring the $(\mathcal{C}^0([-h_M,0];\mathbb{R}^n),\mathcal{H})$-local/global exponential stability property of (\ref{eq: RDE2}) are easily obtained from Theorems~\ref{thm: transfer theorem} and~\ref{thm: final theorem} by assuming the existence of constants $k_1 , k_2 , k_3 > 0$ such that (\ref{eq: main thm 2 - assumed inequality 1}) holds with $\omega_1 (s) = k_1 s^2$ and $\omega_2 (s) = k_2 s^2$ while replacing (\ref{eq: main thm 2 - assumed inequality 3}) by the stronger condition $D^+V_\zeta(t) \leq - k_3 V_\zeta(t)$. Unfortunately, this latter condition is in general difficult to check in practice. Since point-wise dissipation inequalities of the type (\ref{eq: main thm 2 - assumed inequality 3}) are much easier to deal with, an open question is under which extra assumptions such estimates may possibly induce exponential stability properties. In the context of Lyapunov-Karsovskii functionals that solely depend on the system trajectory (i.e., does not depend on its time derivative) and for system trajectories evaluated only in uniform norm, such a question was recently addressed in~\cite{chaillet2019relaxed} in the case of an autonomous and globally Lipschitz continuous vector field.
\end{remark}

We conclude this section with an example of application of Theorem~\ref{thm: final theorem}. This example extends Theorem~\ref{thm: stab sys1} to the case of a vector field presenting a nonlinearity.

\begin{example}
For some constant $h_M > 0$, we study the stability properties of the origin of the following time-delay system:
\begin{subequations}\label{eq: sys2}
\begin{align}
\dot{x}(t) & = M x(t) + N x(t-h(t)) + g(x(t)) , \qquad t \geq 0 \label{eq: sys2 - eq} \\
x(\tau) & = x_0(\tau) ,  \qquad -h_M \leq \tau \leq 0  \label{eq: sys2 - IC}
\end{align}
\end{subequations}
where $M,N \in \mathbb{R}^{n \times n}$, $0 \leq h(t) \leq h_M$ is a continuous time-varying delay, and $x_0 \in \mathcal{C}^0([-h_M,0];\mathbb{R}^n)$ is the initial condition. We assume that $g \in \mathcal{C}^1(\mathbb{R}^n;\mathbb{R}^n)$ with $g(0)=0$ and such that there exist constants $\alpha,\beta,\gamma>0$ for which
\begin{equation}\label{eq: nonlinearity - constraint}
\Vert x \Vert \leq \alpha \quad \Rightarrow \quad \Vert g(x) \Vert \leq \beta \Vert x \Vert^{1+\gamma} .
\end{equation}
Finally we assume that there exist $P_1 , Q \succ 0$ and $P_2,P_3 \in \mathbb{R}^{n \times n}$ such that LMI (\ref{eq: LMI}) holds with $\Gamma = M+N$. 

Introducing the forward invariant set $\mathcal{H} = \{ h\in\mathcal{C}^0(\mathbb{R}_+;\mathbb{R}) \, : \, 0 \leq h(t) \leq h_M \}$, we show that system (\ref{eq: sys2}) is $(\mathcal{C}^0([-h_M,0];\mathbb{R}^n),\mathcal{H})$-locally uniformly asymptotically stable by applying the result of Theorem~\ref{thm: final theorem}. To do so, we introduce the set $A = [0,h_M]$ and define the vector field $f : A \times \mathcal{C}^0([-h_M,0];\mathbb{R}^n) \rightarrow \mathbb{R}^n$ by
\begin{equation*}
f(\xi,\phi) = M \phi(0) + N \phi(-\xi) + g(\phi(0)) .
\end{equation*}
Then $f(\cdot,0)=0$ and it is easy to check that $f$ is continuous. Moreover, one can see that $f$ is Lipchitz continuous with respect to the second argument, uniformly with respect to the first argument, over $A \times \{ \phi \in \mathcal{C}^0([-h_M,0];\mathbb{R}^n) \,:\, \Vert \phi \Vert_\infty \leq H \}$ for any given $H > 0$. Then (\ref{eq: sys2}) is rewritten under the form (\ref{eq: RDE2}). We define the function functional $V : A \times W \rightarrow\mathbb{R}$ by
\begin{equation*}
V(\xi,\phi)
= \phi(0)^\top P_1 \phi(0)
+ \int_{-h_M}^{0} \int_\theta^0 \dot{\phi}(s)^\top Q \dot{\phi}(s) \diff s \diff \theta .
\end{equation*}
Note that this functional is the same as the one used in~\cite{fridman2001new} to prove Theorem~\ref{thm: stab sys1} which deals with the case $g=0$. It is straightforward to show that $V$ is continuous and (\ref{eq: main thm 2 - assumed inequality 1}) holds with $\omega_1(s) = \lambda_m(P_1)s^2$ and $\omega_2(s) = \max(\lambda_M(P_1),h_M \lambda_M(Q))s^2$ where $\lambda_m(S) > 0$ and $\lambda_M(S) > 0$ denote the smallest and largest eigenvalues of the definite positive matrix $S$, respectively.

To conclude, it remains to study the Dini upper right derivative of $V_h(t) = V(h(t),x_t)$ for initial conditions $x_0 \in W$, $h \in \mathcal{H}$, and as soon as $x_t \in Q_H$ for some $H > 0$ to be determined. Noting that
\begin{equation*}
V_h(t) = x(t)^\top P_1 x(t) + \int_{-h_M}^0 \int_{t+\theta}^{t} \dot{x}(s)^\top Q \dot{x}(s) \,\mathrm{d}s\,\mathrm{d}\theta ,   
\end{equation*}
it is easy to show that $V_h$ is continuously differentiable and
\begin{equation*}
\dot{V_h}(t) = 2 x(t)^\top P_1 \dot{x}(t) + h_M \dot{x}(t)^\top Q \dot{x}(t) -\int_{t-h_M}^{t} \dot{x}(s)^\top Q \dot{x}(s) \diff s .
\end{equation*}
Proceeding as in~\cite{fridman2001new} while accounting for the extra term $g(x(t))$, we obtain that
\begin{equation*}
\dot{V}_h(t) \leq \begin{bmatrix} x(t) \\ \dot{x}(t) \end{bmatrix}^\top \Theta \begin{bmatrix} x(t) \\ \dot{x}(t) \end{bmatrix} + 2 \begin{bmatrix} x(t) \\ \dot{x}(t) \end{bmatrix}^\top P^\top \begin{bmatrix} 0 \\ g(x(t)) \end{bmatrix}
\end{equation*}
with $P = \begin{bmatrix} P_1 & 0 \\ P_2 & P_3 \end{bmatrix}$ and
\begin{align*}
\Theta
& =
P^\top \begin{bmatrix} 0 & I \\ \Gamma & -I \end{bmatrix}
+ \begin{bmatrix} 0 & I \\ \Gamma & -I \end{bmatrix}^\top P
+ \begin{bmatrix} 0 & 0 \\ 0 & h_M Q \end{bmatrix} \\
& \phantom{=}\; + h_M P^\top \begin{bmatrix} 0 \\ N \end{bmatrix} Q^{-1} \begin{bmatrix} 0 \\ N \end{bmatrix}^\top P .
\end{align*}
For any $\epsilon > 0$, to be specified later, the use of Young's inequality shows that 
\begin{align*}
\begin{bmatrix} x(t) \\ \dot{x}(t) \end{bmatrix}^\top P^\top \begin{bmatrix} 0 \\ g(x(t)) \end{bmatrix}
& = x(t)^\top P_2^\top g(x(t)) + \dot{x}(t)^\top P_3^\top g(x(t)) \\
& \leq \dfrac{\epsilon}{2} \left( \Vert x(t) \Vert^2 + \Vert \dot{x}(t) \Vert^2 \right)
+ \dfrac{\delta}{2} \Vert g(x(t)) \Vert^2
\end{align*}
with $\delta = \delta(\epsilon) = (\Vert P_2^\top \Vert^2 + \Vert P_3^\top \Vert^2)/\epsilon$ .
Now, based on LMI (\ref{eq: LMI}), the application of the Schur complement yields that $\Theta \prec 0$. Thus we can select the above $\epsilon > 0$ and some $\eta > 0$ such that $\Theta + \epsilon I \preceq -\eta I$. This shows that 
\begin{align*}
\dot{V}_h(t) 
& \leq \begin{bmatrix} x(t) \\ \dot{x}(t) \end{bmatrix}^\top (\Theta+\epsilon I) \begin{bmatrix} x(t) \\ \dot{x}(t) \end{bmatrix} + \delta \Vert g(x(t)) \Vert^2 \\
& \leq -\eta \Vert x(t) \Vert^2 + \delta \Vert g(x(t)) \Vert^2 .
\end{align*}
We introduce the constant $H = \min\left( \alpha , \left( \frac{\eta}{2 \delta \beta^2} \right)^{\frac{1}{2\gamma}} \right) > 0$. Assuming that $x_t \in Q_H$ for a given time $t \geq 0$, we have that $\Vert x(t) \Vert \leq \Vert x_t \Vert_W \leq H$. Since $H \leq \alpha$, we have from (\ref{eq: nonlinearity - constraint}) that $\Vert g(x(t)) \Vert \leq \beta \Vert x(t) \Vert^{1+\gamma}$, while $\Vert x(t) \Vert \leq H \leq \left( \frac{\eta}{2 \delta \beta^2} \right)^{\frac{1}{2\gamma}}$ implies that $\frac{1}{2} \leq 1 - \frac{\delta\beta^2}{\eta} \Vert x(t) \Vert^{2\gamma}$. Consequently, we infer that $x_t \in Q_H$ for a given time $t \geq 0$ implies the following estimates: 
\begin{align*}
\dot{V}_h(t)
& \leq -\eta \Vert x(t) \Vert^2 + \delta \beta^2 \Vert x(t) \Vert^{2(1+\gamma)} \\
& \leq -\eta \Vert x(t) \Vert^2 \left( 1 - \dfrac{\delta \beta^2}{\eta} \Vert x(t) \Vert^{2\gamma} \right) 
\leq -\dfrac{\eta}{2} \Vert x(t) \Vert^2 .
\end{align*}
Then, (\ref{eq: main thm 2 - assumed inequality 3}) holds with $\omega_3(s) = \eta s^2 / 2$. The application of Theorem~\ref{thm: final theorem} shows that system (\ref{eq: sys2}) is $(\mathcal{C}^0([-h_M,0];\mathbb{R}^n),\mathcal{H})$-locally uniformly asymptotically stable.
\end{example}

\section{Conclusion}\label{sec: conclusion}
In the context of retarded differential equations, this note investigated the nature of the stability results obtained via the use of Lyapunov-Krasovskii functionals that depend on the time-derivative of the system trajectory. While the corresponding stability results are generally obtained for absolutely continuous initial conditions with a square-integrable weak derivative and evaluated based on the magnitude of their time-derivative, we showed for a large class of systems the existence of a smoothing effect that allows the extension of these stability properties to continuous initial conditions evaluated in uniform norm. The main feature of the reported result relies on the fact that, for the considered classes of systems, the stability properties of the retarded differential equation, expressed for continuous initial conditions that are evaluated in uniform norm, can be deduced from the study of the system trajectories associated with initial conditions presenting additional regularity assumptions. Since Lyapunov-Krasovskii functionals that depend on the time-derivative of the system trajectory have also proven their efficiency in the assessment of input-to-state stability (ISS) properties for retarded differential equations~\cite{fridman2008input}, a future research direction could include the extension of the results presented in this paper to the ISS framework.



\section*{Acknowledgment}
The authors thank the anonymous reviewers for their comments and suggestions that have led to the nonlinear version of Theorem~\ref{thm: transfer theorem} presented in this paper.

\ifCLASSOPTIONcaptionsoff
  \newpage
\fi



\bibliographystyle{IEEEtranS}
\nocite{*}
\bibliography{IEEEabrv,mybibfile}

\end{document}